\documentclass{amsart}
\usepackage{amsfonts,amsmath,amssymb,graphicx,graphics,color,
lscape,comment,amsbsy,mathrsfs}

\numberwithin{equation}{section}

\newtheorem{theorem}{Theorem}
\numberwithin{theorem}{section}
\newtheorem{lemma}{Lemma}
\numberwithin{lemma}{section}
\newtheorem{prop}{Proposition}
\numberwithin{prop}{section}
\newtheorem{corol}{Corollary}
\numberwithin{corol}{section}
\newtheorem{remark}{Remark}
\numberwithin{remark}{section}

\numberwithin{defi}{section}

\numberwithin{exe}{section}

\newcommand{\notthis}[1]{}

\title[]{\textbf{Inversion of a Class of Singular Integral Operators on Entire Functions}}
\author{R. Nasri(*), A. Simonian(*) and F. Guillemin (**)}
\address{Address:  
(*) Orange Labs, OLN/NMP, Orange Gardens, 
44 avenue de la République, CS 50010, 92326 Chatillon Cedex, France
France (**) Orange Labs Networks Lannion, 2 avenue Pierre Marzin, 22307 Lannion Cedex, Lannion, France}
\email{[ridha.nasri, alain.simonian, fabrice.guillemin]@orange.com}

\begin{document}

\date{Version of \today}

\begin{abstract}

Given constants $x, \nu \in \mathbb{C}$ and the space $\mathscr{H}_0$ of entire functions in 
$\mathbb{C}$ vanishing at $0$, we consider the integro-differential operator 
$$
\mathfrak{L} = \left ( \frac{x \, \nu(1-\nu)}{1-x} \right ) 
\; \delta \circ \mathfrak{M}\, ,
$$
with $\delta = z \, \mathrm{d}/\mathrm{d}z$ and $\mathfrak{M}:\mathscr{H}_0 \rightarrow \mathscr{H}_0$ defined by
$$
\mathfrak{M}f(z) = \int_0^1 e^{-z t^{-\nu}(1-(1-x)t)} \, 
f \left (z \, t^{-\nu}(1-t) \right ) \, \frac{\mathrm{d}t}{t}, 
\qquad z \in \mathbb{C},
$$
for any $f \in \mathscr{H}_0$. Operator $\mathfrak{L}$ originates from an inversion problem in Queuing Theory. Bringing the inversion of $\mathfrak{L}$ back to that of $\mathfrak{M}$ translates into a singular Volterra integral equation, but with no explicit kernel. 

In this paper, the inverse of operator $\mathfrak{L}$ is derived through a new inversion formula recently obtained for infinite matrices with entries involving Hypergeometric polynomials. For 
$x \notin \mathbb{R}^- \cup \{1\}$ and $\mathrm{Re}(\nu) < 0$, we then show that the inverse 
$\mathfrak{L}^{-1}$ of $\mathfrak{L}$ on $\mathscr{H}_0$ has the integral representation 
$$
\mathfrak{L}^{-1}g(z) = \frac{1-x}{2i\pi x} \, e^{z} \int_1^{(0+)} \frac{e^{-xtz}}{t(t-1)} \, 
g \left (z \, (-t)^{\nu}(1-t)^{1-\nu} \right ) \, \mathrm{d}t, 
\qquad z \in \mathbb{C},
$$
for any $g \in \mathscr{H}_0$, where the bounded integration contour in the complex plane starts at 
point 1 and encircles the point 0 in the positive sense. Other related integral representations of 
$\mathfrak{L}^{-1}$ are also provided.

\end{abstract}

\maketitle 


\section{Introduction}


The inversion of an integro-differential operator acting on entire functions in $\mathbb{C}$ is related to a new class of linear inversion formulas with coefficients involving Hypergeometric polynomials. After an overview of the state-of-the-art in the associated fields, we summarize our main contributions.

\subsection{Motivation}
\label{IM}
Consider the following problem:

\textit{\textbf{let constants $x \in \; ]0,1[$, 
$\nu < 0$ and the function $\mathfrak{R}$ defined by}}
\begin{equation}
\mathfrak{R}(\zeta) = x\,\left ( 1 - \zeta \right )^{-\nu}
\left ( 1 - (1-x) \, \zeta \right )^{\nu - 1}, 
\qquad \zeta \in [0,1].
\label{DefRR}
\end{equation}
\textit{\textbf{Let $\mathscr{H}_0$ be the linear space of entire functions in $\mathbb{C}$ vanishing at $z = 0$ and define the integro-differential operator 
$\mathfrak{L}:\mathscr{H}_0 \rightarrow \mathscr{H}_0$ by}}
\begin{equation}
\mathfrak{L}f(z) = \int_0^1 
\biggl[ \left ( 1 + z \mathfrak{R}(\zeta) \right ) 
\, f(\zeta \, \mathfrak{R}(\zeta) \cdot z) - c \, z \, 
\mathfrak{R}(\zeta) \,  
f'(\zeta \, \mathfrak{R}(\zeta) \cdot z) \biggr ] 
e^{- \mathfrak{R}(\zeta) \cdot z} \, \mathrm{d}\zeta 
\label{DefTT}
\end{equation}
\textit{\textbf{for all $z \in \mathbb{C}$, where $f'$ denotes the derivative of $f \in \mathscr{H}_0$ and with the constant $c$ in the integrand equal to}}
$$
c = \frac{1-\nu x}{1-x}.
$$

\textit{\textbf{Given $K \in \mathscr{H}_0$, solve the equation}}
\begin{equation}
\mathfrak{L} \, E^*(z) = K(z), 
\qquad z \in \mathbb{C},
\label{EI0}
\end{equation}
\textit{\textbf{for the unknown $E^* \in \mathscr{H}_0$.}}

This inversion problem has been motivated by an integral equation arising from a problem of Queuing Theory \cite{GQSN18}, namely, the study of the sojourn time in a Processor-Sharing queue with batch customer arrivals.

The operator $\mathfrak{L} = \mathfrak{L}_{x,\nu}$ depends on parameters $x$ and $\nu$. Solving equation (\ref{EI0}) for such parameters is thus equivalent to prove that this operator from 
$\mathscr{H}_0$ to itself is onto. As detailed in this paper, the following \textbf{Properties (I)} and \textbf{(II)} for $\mathfrak{L}$ and the associated equation (\ref{EI0}) can be successively outlined:

\textbf{(I) Reduction to a Linear System: power series expansions}
\begin{equation}
E^*(z) = \sum_{\ell = 1}^{+\infty} E_\ell \, \frac{z^\ell}{\ell!}, \; \; \; \; 
K(z) = \sum_{b = 1}^{+\infty} (-1)^b K_b \, \frac{z^b}{b!}, 
\qquad z \in \mathbb{C},
\label{defE*}
\end{equation}
\textbf{for a solution $E^* \in \mathscr{H}_0$ and the given 
$K \in \mathscr{H}_0$ reduce the resolution of (\ref{EI0}) to that of the infinite lower-triangular linear system}
\begin{equation}
\forall \, b \in \mathbb{N}^*, \qquad 
\sum_{\ell = 1}^b (-1)^\ell \binom{b}{\ell} \, 
Q_{b,\ell} \, E_{\ell} = K_b,  
\label{T0}
\end{equation}
\textbf{with unknown $E_\ell$, $\ell \in \mathbb{N}^*$, and where the coefficient matrix $Q = (Q_{b,\ell})_{b, \ell \in \mathbb{N}^*}$, on account of the specific function $\mathfrak{R}$ introduced in 
(\ref{DefRR}), is given by}
\begin{equation}
Q_{b,\ell} = \frac{\Gamma(b)\Gamma(1-b\nu)}{\Gamma(b-b\nu)} \, 
\frac{x}{x-1} \; F(\ell - b,-b \nu;-b;x), 
\qquad 1 \leqslant \ell \leqslant b.
\label{Q0}
\end{equation}

In (\ref{Q0}), $\Gamma$ is the Euler Gamma function and 
$F(\alpha,\beta;\gamma;\cdot)$ denotes the Gauss Hypergeometric function with complex parameters $\alpha$, $\beta$, $\gamma \notin -\mathbb{N}$. Recall that $F(\alpha,\beta;\gamma;\cdot)$ reduces to a polynomial with degree $-\alpha$ (resp. $-\beta$) if $\alpha$ (resp. $\beta$) equals a non positive integer; expression (\ref{Q0}) for coefficient $Q_{b,\ell}$ thus involves a Hypergeometric polynomial with degree $b - \ell$ in both arguments $x$ and $\nu$.

The diagonal coefficients $Q_{b,b}$, $b \geqslant 1$, are non-zero so that lower-triangular system (\ref{T0}) has a unique solution;  equivalently, this proves the uniqueness of the solution $E^* \in \mathscr{H}_0$ to (\ref{EI0}). To make this solution explicit in terms of parameters, write system (\ref{T0}) equivalently as 
\begin{equation}
\forall \, b \in \mathbb{N}^*, \qquad 
\sum_{\ell = 1}^b A_{b,\ell}(x,\nu) \, E_\ell = \widetilde{K}_b,
\label{T0BIS}
\end{equation}
with the reduced right-hand side $(\widetilde{K}_b)$ defined by 
$$
\widetilde{K}_b = \frac{\Gamma(b-b\nu)}{\Gamma(b)\Gamma(1-b\nu)} \frac{x-1}{x}\cdot K_b, 
\qquad b \geqslant 1,
$$ 
and with matrix $A(x,\nu) = (A_{b,\ell}(x,\nu))$ given by
\begin{equation}
A_{b,\ell}(x,\nu) = (-1)^\ell \binom{b}{\ell} F(\ell-b,-b\nu;-b;x), 
\qquad 1 \leqslant \ell \leqslant b.
\label{T0TER}
\end{equation}
As recently shown \cite{NAS20}, the linear relation (\ref{T0BIS}) to which initial system (\ref{T0}) has been recast can be explicitly inverted for any right-hand side $(K_b)_{b \in \mathbb{N}^*}$, the inverse matrix $B(x,\nu) = A(x,\nu)^{-1}$ involving also Hypergeometric polynomials as well. This consequently solves system (\ref{T0}) explicitly, hence integral equation (\ref{EI0});

\textbf{(II) Factorization: operator $\mathfrak{L}$ can be factored as}
\begin{equation}
\mathfrak{L} =  \frac{x\nu(1-\nu)}{1-x} \cdot \delta \circ \mathfrak{M}
\label{DefTTbis}
\end{equation}
\textbf{where $\delta = z \, \mathrm{d}/\mathrm{d}z$ and $\mathfrak{M}$ is the integral operator defined by}
\begin{equation}
\mathfrak{M}f(z) = 
\int_0^1 e^{-z t^{-\nu}(1-(1-x)t)} \, 
f \left (z \, t^{-\nu}(1-t) \right ) \, \frac{\mathrm{d}t}{t}, 
\qquad z \in \mathbb{C},
\label{DefSS}
\end{equation}
\textbf{for all $f \in \mathscr{H}_0$.} Using factorization (\ref{DefTTbis}), the resolution of (\ref{EI0}) is thus equivalent to solving equation 
\begin{equation}
\mathfrak{M}E^* = K_1
\label{DefK1}
\end{equation}
with right-hand side
$$
K_1(z) = \frac{1-x}{\nu(1-\nu)x} \cdot \int_0^z \frac{K(\zeta)}{\zeta} \, \mathrm{d}\zeta, 
\qquad z \in \mathbb{C},
$$
where $K_1 \in \mathscr{H}_0$ as soon as $K \in \mathscr{H}_0$. Integral equation (\ref{DefK1}) can be in turn recast into the Volterra equation 
\begin{equation}
\int_0^{\widehat{\tau} z} 
\Psi \left (z, \frac{\xi}{z}\right ) \, E^*(\xi) \, \mathrm{d}\xi = 
z\cdot K_1(z), \qquad z \in \mathbb{C},
\label{Fredh0}
\end{equation}
for some constant $\widehat{\tau}$ and a kernel $\Psi(z,\cdot)$. As $\Psi(z,\cdot)$ has an integrable singularity of order $\Psi(z,\tau) = O \left (\widehat{\tau} - \tau \right )^{-1/2}$ near point 
$\tau = \widehat{\tau}$, (\ref{Fredh0}) is therefore a singular Volterra integral equation of the first kind.

\subsection{State-of-the-art}
As equation (\ref{EI0}) or (\ref{DefK1}) can be recast into the singular Volterra equation of the first kind (\ref{Fredh0}), we here briefly review known results for this class of integral equations. 

Given the constant $\alpha \in \; ]0,1[$, the standard case for such singular equations is that of the classical Abel's equation
$$
\int_0^z \frac{E(\xi)}{(z-\xi)^\alpha} \, \mathrm{d}\xi = \kappa(z), 
\qquad z \in [0,r],
$$
on a real interval $[0,r]$, for the unknown function $E$ and some 
given function $\kappa$ (see \cite[Chap.7]{BIT95}, \cite[Chap.2]
{EST00}, \cite[Chap.1]{GOR91}. If $\kappa$ is absolutely continuous on 
$[0,r]$, then Abel's equation has the unique solution 
$E \in \mathrm{L}^1[0,r]$ given by
\begin{align}
E(z) = & \, \frac{\sin(\pi\alpha)}{\pi} \cdot 
\frac{\mathrm{d}}{\mathrm{d}z} 
\left [ \int_0^z \frac{\kappa(\xi)}{(z-\xi)^{1-\alpha}} \, 
\mathrm{d}\xi \right ]
\nonumber \\
= & \, \frac{\sin(\pi\alpha)}{\pi} 
\left [ \frac{\kappa(0)}{z^{1-\alpha}} + 
\int_0^z \frac{\kappa'(\xi)}{(z-\xi)^{1-\alpha}} \, 
\mathrm{d}\xi \right ], 
\qquad z \in [0,r].
\label{SolAbel}
\end{align}
This solution extends to a complex variable $z \in \mathbb{C}$ pertaining to a neighborhood of point $0$ where function $\kappa$ is assumed to be analyti; the solution $E$ is then analytic in a neighborhood of $z = 0$ if condition $\kappa(0) = 0$ holds, that is, if and only if $\kappa \in \mathscr{H}_0$. 

More generally, let a compact subset $\Omega \subset \mathbb{C}$ and a singular operator 
$\mathfrak{J}:E \mapsto \mathfrak{J}E$ defined by
$$
\mathfrak{J}E(z) = \int_0^z N(z,\xi)E(\xi) \, \mathrm{d}\xi, 
\qquad z \in \Omega,
$$
where the kernel $N$ verifies
$$
\vert N(z,\xi) \vert \leqslant \frac{M}{\vert z - \xi \vert^\alpha}, \qquad 
z, \; \xi \in \Omega, \; z \neq \xi,
$$
for some constant $M > 0$ and $\alpha \in \; ]0,1[$. Operator $\mathfrak{J}$ is known to be continuous (and also compact) on space $\mathscr{C}^0[\Omega]$ \cite[Theorem 2.29]{KRESS14}. No general results are available, however, on the inverse of $\mathfrak{J}$ on a subspace of 
$\mathscr{C}^0[\Omega]$.

This standard framework may nevertheless suggest the existence of an integral representation of the kind (\ref{SolAbel}) for the entire solution to either equation (\ref{EI0}) (\ref{DefK1}) or (\ref{Fredh0}). In this paper, we will show how such integral representations can be obtained for the solutions of these singular equations. 

\subsection{Paper contribution}
The main contributions of this paper can be summarized as follows:

\textbf{(A)} we first prove the above mentioned \textbf{Reduction Property I} (Section  \ref{SecA1}) whereby integral equation (\ref{EI0}) is reduced to linear system (\ref{T0}) with coefficients related to Hypergeometric polynomials; 

\textbf{(B)} we next justify the \textbf{Factorization Property II} (Section \ref{SecA2}) for the integro-differential  operator $\mathfrak{L}$. We further specify how equation (\ref{DefK1}) can be recast into a Volterra integral equation with singular kernel (Section \ref{SecA3});

\textbf{(C)} the previous results finally enable us to derive an integral representation of the inverse 
$\mathfrak{L}^{-1}$ of operator $\mathfrak{L}$ in space $\mathscr{H}_0$ in the form of the contour integral
$$
\mathfrak{L}^{-1}g(z) = \frac{1-x}{2i\pi x} \, e^{z} \int_1^{(0+)} \frac{e^{-xtz}}{t(t-1)} \, 
g \left (z \, (-t)^{\nu}(1-t)^{1-\nu} \right ) \, \mathrm{d}t, 
\qquad z \in \mathbb{C},
$$
for any $g \in \mathscr{H}_0$, where the finite contour in the complex plane starts at 1 and encircles 0 in the positive sense (Section \ref{SecB1}). By using suitable variable change in the latter, other related integral representations of the inverse $\mathfrak{L}^{-1}$ are also provided.


\section{Preliminaries}
\label{LTS}


\subsection{Infinite matrices inversion}
\label{IMI}
We first recall the results established in \cite{NAS20} for the inversion of some class of lower-triangular matrices. These results will be used below for the inversion of operator $\mathfrak{L}$. 

While stated for general matrices $\mathbf{A}(x,\nu; \alpha,\beta,\gamma)$ depending on $x$, $\nu$ and three other complex parameters $\alpha$, $\beta$ and $\gamma$ \cite[Theorem 2.3]{NAS20}, we here only use the inversion property particularized to the matrix $\mathbf{A}(x,\nu)$ introduced in 
(\ref{T0TER}) and corresponding to the sub-case $\alpha = \beta = \gamma = 0$. In such a case, the inversion formula for lower-triangular matrices involving Hypergeometric polynomials 
$F(m,\cdot;\cdot;x)$, $m \in -\mathbb{N}$, can be stated as follows. 

\begin{theorem} (\cite[Sect. 2.2]{NAS20})
\textbf{Let $x, \nu \in \mathbb{C}$ and define the lower-triangular matrices 
$A(x,\nu)$ and $B(x,\nu)$ by}
	\begin{equation}
	\left\{
	\begin{array}{ll}
	A_{n,k}(x,\nu) = \displaystyle (-1)^k\binom{n}{k}F(k-n,-n\nu;-n;x),
  \\ \\
	B_{n,k}(x,\nu) = \displaystyle (-1)^k\binom{n}{k}F(k-n,k\nu;k;x)
	\end{array} \right.
	\label{DefABxNU}
	\end{equation}
\textbf{for $1 \leqslant k \leqslant n$. The inversion formula}
	\begin{equation}
	T_n = \sum_{k=1}^{n}A_{n,k}(x,\nu)S_k 
	\Longleftrightarrow 
	S_n=\sum_{k=1}^{n}B_{n,k}(x,\nu)T_k, \quad n \in \mathbb{N}^*,
	\label{eq:inversionR}
	\end{equation}
\textbf{holds for any pair of complex sequences 
$(S_n)_{n \in \mathbb{N}^*}$ and $(T_n)_{n \in \mathbb{N}^*}$.}
	\label{PropIn}
\end{theorem}

As a direct consequence of Theorem \ref{PropIn}, a remarkable functional identity can be derived for the exponential generating functions of sequences related by the inversion formula.

\begin{corol} (\cite[Sect. 3.2]{NAS20})
\textbf{Given sequences $S$ and $T$ related by the inversion formulae 
$S = B(x,\nu) \cdot T \Leftrightarrow T = A(x,\nu) \cdot S$, the exponential generating function $\mathfrak{G}_S^*$ of the sequence $S$ can be expressed by}
\begin{equation}
\mathfrak{G}_S^*(z) = \exp(z) \cdot 
\sum_{k \geqslant 1} (-1)^k T_k \, \frac{z^k}{k!} \, \Phi(k\nu;k;-x \, z), 
\qquad z \in \mathbb{C},
\label{GsExp}
\end{equation}
\textbf{where $\Phi(\alpha;\beta;\cdot)$ denotes the Confluent Hypergeometric function with parameters $\alpha$, 
$\beta \notin -\mathbb{N}$.}
\label{corEGF}
\end{corol}

\subsection{Parameters range}
Operator $\mathfrak{L}$ has been initially introduced for real parameters $x \in \; ]0,1[$ and $\nu < 0$. Following the results recalled in Section \ref{IMI} and stated for arbitrary complex parameters, we hereafter extend the definition (\ref{DefTT}) of $\mathfrak{L}$ to complex values, namely  
\begin{itemize}
\item $x \in \; \mathbb{C} \setminus 
( \mathbb{R}^- \cup \{1\})$ (so that $1/(1-x)$ is finite and does not belong to the integration interval $[0,1]$) 
\item and $\nu \in \mathbb{C}$ such that $\mathrm{Re}(\nu) < 0$. 
\end{itemize}
Within these assumptions, it is easily verified that $\mathfrak{L}(\mathscr{H}_0) \subset \mathscr{H}_0$ where $\mathscr{H}_0$ is again the linear space of entire functions in 
$\mathbb{C}$ vanishing at $0$.

\begin{remark}
Operator $\mathfrak{L}$, well-defined for $\mathrm{Re}(\nu) < 0$, may not exist for other values of 
$\nu$. In fact, consider the function $f \in \mathscr{H}_0$ defined by 
$f(z) = z \, e^{(1-x)z}$, $z \in \mathbb{C}$. By definition (\ref{DefTT}), it is easily verified that, for 
$\mathrm{Re}(\nu) < 0$, its image $\mathfrak{L}f \in \mathscr{H}_0$ is given by
$$
\mathfrak{L}f(z) = - \, \frac{z}{1-x} \, 
\Phi \left ( 1 - \frac{1}{\nu};1 - \frac{1}{\nu}; -z \right ), 
\qquad z \in \mathbb{C},
$$
$\Phi(\cdot;\cdot;\cdot)$ denoting the Kummer Confluent Hypergeometric function. For 
$\mathrm{Re}(\nu) > 0$, however, its image is given by
$$
\mathfrak{L}f(z) = - \, \frac{1-\nu}{\nu(1-x)} \, 
z^{\frac{1}{\nu}} \, \Gamma \left ( 1 - \frac{1}{\nu}; z \right ), 
\qquad z \neq 0
$$
(where $\Gamma(\cdot;\cdot)$ is the incomplete Gamma function), so that 
$\mathfrak{L}f \notin \mathscr{H}_0$ in this case.
\end{remark}


\section{Properties of operator $\mathfrak{L}$}
\label{SecA}


\subsection{Reduction to a linear system}
\label{SecA1}
We have claimed in \ref{IM}.\textbf{(I)} that the integro-differential equation (\ref{EI0}) reduces to the infinite system (\ref{T0}). We justify this assertion by showing how the coefficients of system (\ref{T0}) can be expressed in terms of Hypergeometric polynomials.

\begin{prop}
\textbf{The Reduction Property (I) holds, that is, equation (\ref{EI0}) reduces to system (\ref{T0}) with matrix $Q = (Q_{b,\ell})_{1 \leqslant \ell  \leqslant b}$ related to Hypergeometric polynomials as given in (\ref{Q0}).}
\label{LemT0}
\end{prop}

\begin{proof}
To derive system (\ref{T0}), we expand both sides of (\ref{EI0}) into power series of variable $z$ and identify like powers on each side. The series expansion (\ref{defE*}) of $E^*(Z)$ in powers of $Z$ first provides
\begin{equation}
(1 + Z)E^*(\zeta Z) - c \, Z \, \frac{\mathrm{d} E^*}{\mathrm{d} z}(\zeta Z) = 
\sum_{b \geqslant 1} \Lambda_b(\zeta) \frac{{Z}^b}{b!}
\label{Bsv00}
\end{equation}
where we set
$\Lambda_b(\zeta) = \zeta^b E_b + b \, \zeta^{b-1} 
E_{b-1} - b \, c \zeta^{b-1}E_b$ 
for all $\zeta$ and with the constant $c = (1-\nu x)/(1-x)$; applying equality (\ref{Bsv00}) to the argument $Z = \mathfrak{R}(\zeta) \cdot z$, the integrand of $\mathfrak{L} E^*(z)$ in (\ref{DefTT}) can then be expanded into a power series of $z$ as
\begin{equation}
\mathfrak{L} E^*(z) = \int_0^{1} \biggl[ \sum_{b \geqslant 1} 
\Lambda_b(\zeta) \, \frac{\zeta^b \, \mathfrak{R}(\zeta)^b \, z^b}{b!} \biggr]
\, e^{-\mathfrak{R}(\zeta) \, z} \, \mathrm{d}\zeta.
\label{Bsv0}
\end{equation}
Now, expanding the exponential $e^{-\mathfrak{R}(\zeta) \, z}$ of the integrand in (\ref{Bsv0}) into a power series of $z$ gives the expansion
\begin{equation}
\mathfrak{L} E^*(z) = \sum_{b \geqslant 0} (-1)^b \frac{z^b}{b!} \, 
\sum_{\ell=0}^b (-1)^{\ell} \binom{b}{\ell} \int_0^{1} 
\zeta^\ell \Lambda_\ell(\zeta) \mathfrak{R}(\zeta)^b \, \mathrm{d}\zeta
\label{Bsv}
\end{equation}
(after noting that $\Lambda_0(\zeta) = 0$ since $E_0 = 0$ by definition). On account of expansion (\ref{Bsv}) with the above definition (\ref{Bsv00}) of $\Lambda_\ell(\zeta)$, together with the expansion (\ref{defE*}) for $K(z)$, the identification of like powers of these expansions readily yields the relation
\begin{equation}
\sum_{\ell=1}^b (-1)^\ell \binom{b}{\ell} B_{b,\ell} E_\ell^* + 
\sum_{\ell=1}^b (-1)^\ell \binom{b}{\ell} \ell \, M_{b,\ell-1} 
E_{\ell-1} = K_b, \qquad b \geqslant 1,
\label{T0bis}
\end{equation}
with $B_{b,\ell} = M_{b,\ell} - \ell \, c \, M_{b,\ell-1}$, 
where $M_{b,\ell}$ denotes the definite integral
\begin{equation}
M_{b,\ell} = \int_0^1 \zeta^\ell \, \mathfrak{R}(\zeta)^b \, \mathrm{d}\zeta, 
\qquad 1 \leqslant \ell \leqslant b.
\label{DefMbell}
\end{equation}
By first changing the index in the second sum in the left-hand side of 
(\ref{T0bis}) and then using identity 
$\binom{b}{\ell+1} = (b-\ell) \cdot \binom{b}{\ell}/(\ell + 1)$, 
(\ref{T0bis}) reduces to (\ref{T0}) with coefficients
\begin{equation}
Q_{b,\ell} = (\ell + 1 - b) M_{b,\ell} - 
\ell \, c \, M_{b,\ell-1}, \qquad 1 \leqslant \ell \leqslant b.
\label{DefQbell}
\end{equation}
The calculation of integral $M_{b,\ell}$ in (\ref{DefMbell}) in terms of Hypergeometric functions and its reduction to Hypergeometric polynomials is detailed in Appendix \ref{A4}; this eventually provides expression 
(\ref{Q0}) for the coefficients of matrix $Q = (Q_{b,\ell})$. 
\end{proof}

We can now deduce the unique solution to system (\ref{T0}).

\begin{corol}
\textbf{Let $\nu \in \mathbb{C}$ with $\mathrm{Re}(\nu) < 0$. Given the sequence 
$(K_b)_{b \geqslant 1}$, the unique solution 
$(E_b)_{b \geqslant 1}$ to system (\ref{T0}) is given by}
\begin{equation}
E_b = \frac{x-1}{x} \sum_{\ell = 1}^b (-1)^{\ell } 
\binom{b}{\ell} F(\ell-b,\ell \nu;\ell;x) \, 
\frac{\Gamma(\ell - \ell \nu)}{\Gamma(\ell)\Gamma(1-\ell \nu)} \, K_\ell
\label{E0}
\end{equation}
\textbf{for all $b \geqslant 1$.}
\label{C2}
\end{corol}

\begin{proof}
By expression (\ref{Q0}) for the coefficients of lower-triangular matrix $Q$, equation (\ref{T0}) equivalently reads
\begin{equation}
\sum_{\ell=1}^b (-1)^\ell \binom{b}{\ell} F(\ell-b,-b\nu;-b;x) \cdot 
E_\ell = \widetilde{K}_b, \qquad 1 \leqslant \ell \leqslant b,
\label{S0}
\end{equation}
when setting
\begin{equation}
\widetilde{K}_b = \displaystyle \, \frac{\Gamma(b - b \nu)}
{\Gamma(b)\Gamma(1-b \nu)} \frac{x-1}{x} \cdot K_b, 
\qquad b \geqslant 1.
\label{S1}
\end{equation}
The application of inversion Theorem \ref{PropIn} to lower-triangular system (\ref{S0}) readily provides the solution sequence 
$(E_\ell)_{\ell \in \mathbb{N}}$ in terms 
of the sequence $(\widetilde{K}_b)_{b \in \mathbb{N}^*}$; using then  transformation (\ref{S1}), the final solution (\ref{E0}) for the sequence 
$(E_\ell)_{\ell \in \mathbb{N}^*}$ follows. 
\end{proof}

\subsection{Factorization of $\mathfrak{L}$}
\label{SecA2}
We now prove the factorization property \ref{IM}.\textbf{(II)} for integro-differential operator 
$\mathfrak{L}$.

\begin{prop}
\textbf{The Factorization Property (II) holds, that is, the linear operator $\mathfrak{L}$ on space 
$\mathscr{H}_0$ can be factored as in (\ref{DefTTbis}) in terms of operators 
$\delta = z \, \mathrm{d}/\mathrm{d}z$ and $\mathfrak{M}$.}
\label{FactorP}
\end{prop}

\begin{proof}
Calculating the exponential generating function of the sequence $(-1)^bK_b$, $b \geqslant 1$, from relation (\ref{T0}) with help of (\ref{Q0}) for the coefficients of matrix $Q$ gives
\begin{align}
& \mathfrak{L}E^*(z) = K(z) = \sum_{b \geqslant 1} (-1)^b K_b \, z^b \; = 
\nonumber \\
& \sum_{b \geqslant 1} - \frac{\Gamma(b)\Gamma(1-b\nu)}{\Gamma(b-b\nu)} \, \frac{x}{1-x} \, \frac{(-z)^b}{b!} \sum_{\ell = 1}^b (-1)^\ell \binom{b}{\ell} \, F(\ell-b,-b\nu;-b;x) E_\ell,
\nonumber
\end{align}
for all $z \in \mathbb{C}$, that is,
\begin{align}
\mathfrak{L}E^*(z) = 
& \sum_{\ell \geqslant 1} \frac{(-1)^\ell}{\ell!} E_\ell \; \times 
\nonumber \\
& \sum_{b \geqslant \ell} - \frac{\Gamma(b)\Gamma(1-b\nu)}{\Gamma(b-b\nu)} \, 
\frac{x}{1-x} \, \frac{(-z)^b}{(b-\ell)!} F(\ell-b,-b\nu;-b;x)
\label{F0}
\end{align}
(after changing the summation order on indexes $b$ and $\ell$). Applying the general identity (\ref{IDFG}) to parameters $m = b-\ell \geqslant 0$, 
$\beta = -b\nu$ and $\gamma = \ell - b\nu +1$ to express polynomial 
$F(\ell-b,-b\nu,-b;x)$ in terms of polynomial 
$F(\ell-b,-b\nu,\ell-b\nu+1;1-x)$, we further obtain
\begin{align}
& - \frac{\Gamma(b)\Gamma(1-b\nu)}{\Gamma(b - b\nu)} \, F(\ell-b,-b\nu,-b;x) 
\; = 
\nonumber \\
& - \frac{(1-\nu)\Gamma(\ell+1)\Gamma(1-b\nu)}{\Gamma(\ell-b\nu+1)} \, 
F(\ell-b,-b\nu,\ell-b\nu+1;1-x);
\label{F1}
\end{align}
using the integral representation recalled in Appendix \ref{A4} - 
Equ.(\ref{HyperGauss}) for the factor $F(\ell-b,-b\nu,\ell-b\nu+1;1-x)$ in the right-hand side of (\ref{F1}) eventually yields 
$$
- \frac{\Gamma(b)\Gamma(1-b\nu)}{\Gamma(b - b\nu)} \, F(\ell-b,-b\nu,-b;x) 
= b\nu(1-\nu) \int_0^1 t^{-b\nu-1}(1-t)^\ell(1-(1-x)t)^{b-\ell} \, 
\mathrm{d}t.
$$
Now, replacing the latter into the right-hand side of (\ref{F0}) provides
\begin{equation}
\mathfrak{L}E^*(z) = \frac{x\nu(1-\nu)}{1-x} \, 
\sum_{\ell \geqslant 1} 
\frac{(-1)^\ell}{\ell!} \, E_\ell \cdot S_{\ell}(z)
\label{F3}
\end{equation}
where
$$
S_{\ell}(z) = \sum_{b \geqslant \ell} \frac{b}{(b-\ell)!} \cdot 
\left (-z\right )^b \, 
\int_0^1 t^{-b\nu-1}(1-t)^\ell(1-(1-x)t)^{b-\ell} \, \mathrm{d}t.
$$
With the index change $b' = b - \ell$, the sum $S_{\ell}(z)$ for given 
$\ell$ equivalently reads
\begin{align}
S_\ell(z) = & \, \sum_{b' \geqslant 0} \frac{b'+\ell}{b'!} \cdot 
\left (-z\right )^{b'+\ell} \, 
\int_0^1 t^{-(b'+\ell)\nu}(1-t)^\ell(1-(1-x)t)^{b'} \, 
\frac{\mathrm{d}t}{t} 
\nonumber \\
= & \, z \, \frac{\mathrm{d}}{\mathrm{d}z} \left [ 
\sum_{b' \geqslant 0} \frac{1}{b'!} \cdot 
\left (-z\right )^{b'+\ell} \, 
\int_0^1 t^{-(b'+\ell)\nu}(1-t)^\ell(1-(1-x)t)^{b'} \, 
\frac{\mathrm{d}t}{t} \right ] 
\nonumber
\end{align}
hence
\begin{equation}
S_\ell(z) = z \, \frac{\mathrm{d}}{\mathrm{d}z} 
\left [ \left ( -z \right )^{\ell} 
\int_0^1 t^{-\ell \nu} (1-t)^\ell \, \frac{\mathrm{d}t}{t} \times 
\exp \left ( - z \, t^{-\nu}[1-(1-x)t] \right ) \right ].
\label{F4}
\end{equation}
Replacing expression (\ref{F4}) into the left-hand side of (\ref{F3}), the linearity of operator $\delta = z \, \mathrm{d}/\mathrm{d}z$ and the permutation of the summation on index $\ell$ with the integration with respect to variable $t \in [0,1]$ enable us to obtain
$$
\mathfrak{L}E^*(z) = \frac{x\nu(1-\nu)}{1-x} \cdot  
\delta \Biggl [ \int_0^1 \frac{\mathrm{d}t}{t} \, 
e^{-z\, t^{-\nu}[1-(1-x)t]} \; 
\sum_{\ell \geqslant 1} 
\frac{E_\ell}{\ell!} 
z^{\ell} t^{-\ell \nu}(1-t)^\ell \Biggr ],
$$
that is,
\begin{align}
\mathfrak{L}E^*(z) = & \, \frac{x\nu(1-\nu)}{1-x} \cdot  
\delta \Biggl [ \int_0^1 \frac{\mathrm{d}t}{t} 
e^{-z\, t^{-\nu}[1-(1-x)t]} \times 
E^* \left (z\, t^{-\nu}(1-t) \right ) \Biggr ] 
\nonumber \\
= & \, \frac{x\nu(1-\nu)}{1-x} \cdot (\delta \circ \mathfrak{M})E^*(z), 
\qquad z \in \mathbb{C},
\nonumber
\end{align}
for any function $E^* \in \mathscr{H}_0$, as claimed in (\ref{DefTTbis}), with the corresponding definition of integral operator $\mathfrak{M}$ on space 
$\mathscr{H}_0$. 
\end{proof}

\subsection{The Volterra equation}
\label{SecA3}
\begin{figure}[htb]
\scalebox{0.9}{\includegraphics[width=13cm, 
trim = 0cm 10.5cm 11cm 1.5cm,clip]
{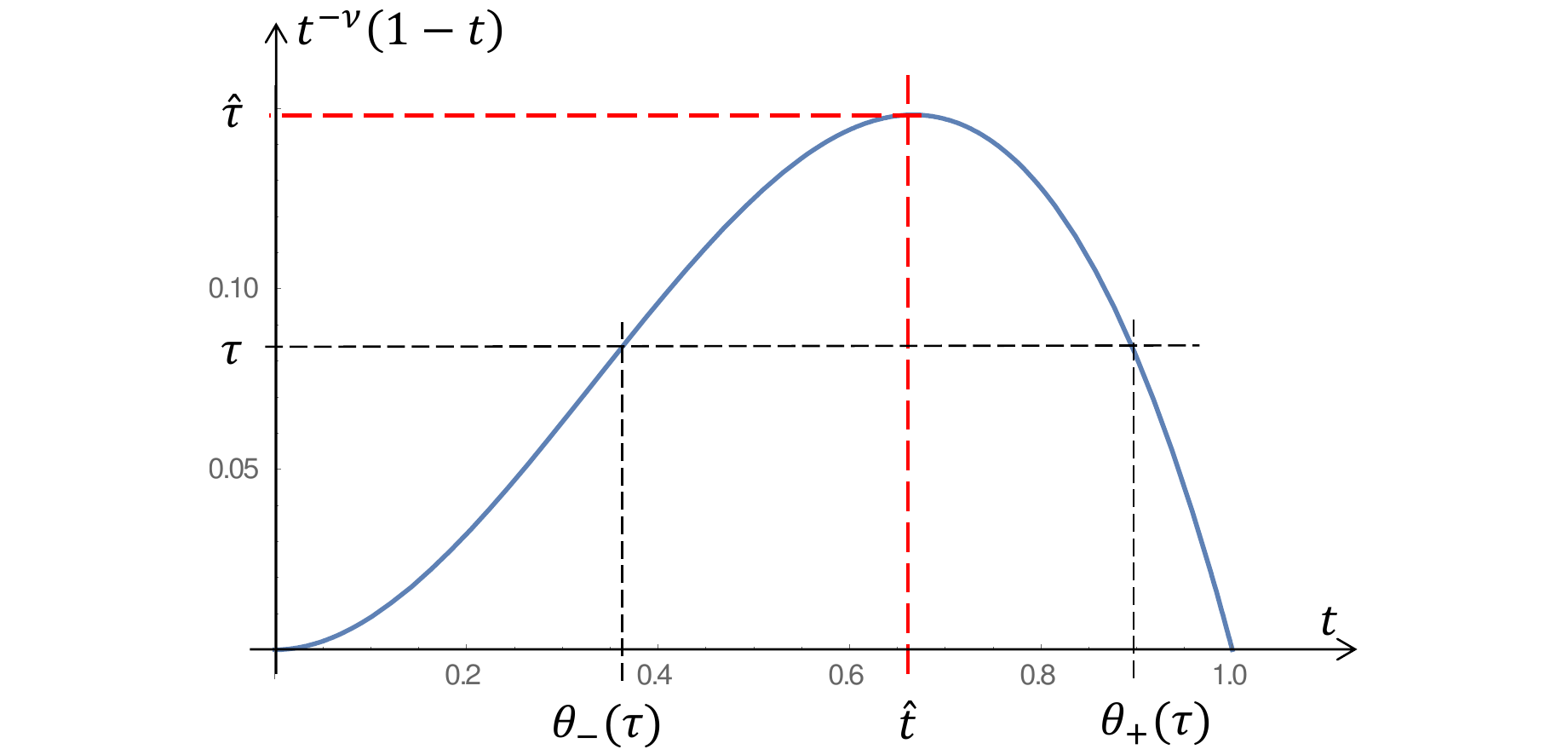}}
\caption{\textit{Graph of function $\tau:t \in [0,1] \mapsto t^{-\nu}(1-t)$ 
(here $\nu = -2$ for illustration).}}
\label{GraphTau}
\end{figure}

As outlined in the Introduction, the factorization (\ref{DefTTbis}) of operator $\mathfrak{L}$ allows one to write equation (\ref{EI0}) equivalently as 
\begin{equation}
\mathfrak{M}E^* = K_1
\label{EI0bis}
\end{equation}
where $K_1 \in \mathscr{H}_0$ relates to the initial function $K$ as in (\ref{DefK1}). As the real function
$\tau:t \in [0,1] \mapsto t^{-\nu}(1-t)$ has a unique maximum at point $\widehat{t} = \nu/(\nu - 1)$ for $\nu < 0$, we can introduce the variable changes 
$t \in [0,\widehat{t}] \mapsto \tau_-(t) = t^{-\nu}(1-t)$ and 
$t \in [\widehat{t},1] \mapsto \tau_+(t) = t^{-\nu}(1-t)$ on segments $[0,\widehat{t}]$ and 
$[\widehat{t},1]$, respectively; we further denote by
\begin{equation}
\theta_- = \tau_-^{-1}, \qquad \theta_+ = \tau_+^{-1}
\label{Thetapm}
\end{equation}
the respective inverse mappings of $\tau_-$ and $\tau_+$, both defined on segment $[0,\widehat{\tau}]$ where $\widehat{\tau} = \tau(\widehat{t}) = (\widehat{t})^{-\nu}(1-\widehat{t})$ (see illustration on Fig.\ref{GraphTau}). The variable changes $\tau_-$ and $\tau_+$ then allow us to write equation (\ref{EI0bis}) as a singular Volterra equation.

\begin{corol}
\textbf{Given constants $x$ and $\nu$ as above, the equivalent equation 
(\ref{EI0bis}) can be recast into the singular Volterra integral equation}
\begin{equation}
\int_0^{\widehat{\tau}\, z} 
\left [ \Psi_- \left ( z, \frac{\xi}{z}\right ) - 
\Psi_+ \left ( z, \frac{\xi}{z} \right ) \right ] 
E^*(\xi) \, \mathrm{d}\xi = 
z\cdot K_1(z), \qquad z \in \mathbb{C},
\label{VOL0}
\end{equation}
\textbf{where we set}
$$
\Psi_{\pm}(z,\tau) = 
\frac{e^{-z
\theta_{\pm}(\tau)^{-\nu}(1-(1-x)\theta_{\pm}(\tau))}}
{\theta_{\pm}(\tau)^{-\nu}(-\nu+(\nu-1)\theta_{\pm}(\tau))}, 
\qquad 0 \leqslant \tau \leqslant \widehat{\tau},
$$
\textbf{with $\theta_{\pm}$ introduced in (\ref{Thetapm}), and where 
$K_1 \in \mathscr{H}_0$ is defined by (\ref{DefK1}).}
\label{C3}
\end{corol}

\noindent
We refer to Appendix \ref{A5} for the proof of Corollary \ref{C3}. It is noted there that the kernel 
$\tau \mapsto \Psi_-(z,\tau) - \Psi_+(z,\tau)$ of Volterra equation (\ref{VOL0}) is singular with an integrable singularity at the boundary $\tau = \widehat{\tau}$ of order 
$O(\widehat{\tau} - \tau)^{-1/2}$. 

Although giving a reformulation to initial equation (\ref{EI0}), equation (\ref{VOL0}) remains difficult to solve as its kernel depends on inverse functions $\theta_-$ and $\theta_+$ which cannot be made explicit simply.


\section{Inversion of operator $\mathfrak{L}$}
\label{SecB1}


We now provide integral representations for the inverse of operator 
$\mathfrak{L}$ on $\mathscr{H}_0$ or, equivalently, integral representations for the solution of integral equation (\ref{EI0}) addressed in the Introduction.

\begin{figure}[htb]
\scalebox{0.8}{\includegraphics[width=14cm, trim = 4cm 7.5cm 3cm 
4.2cm,clip]
{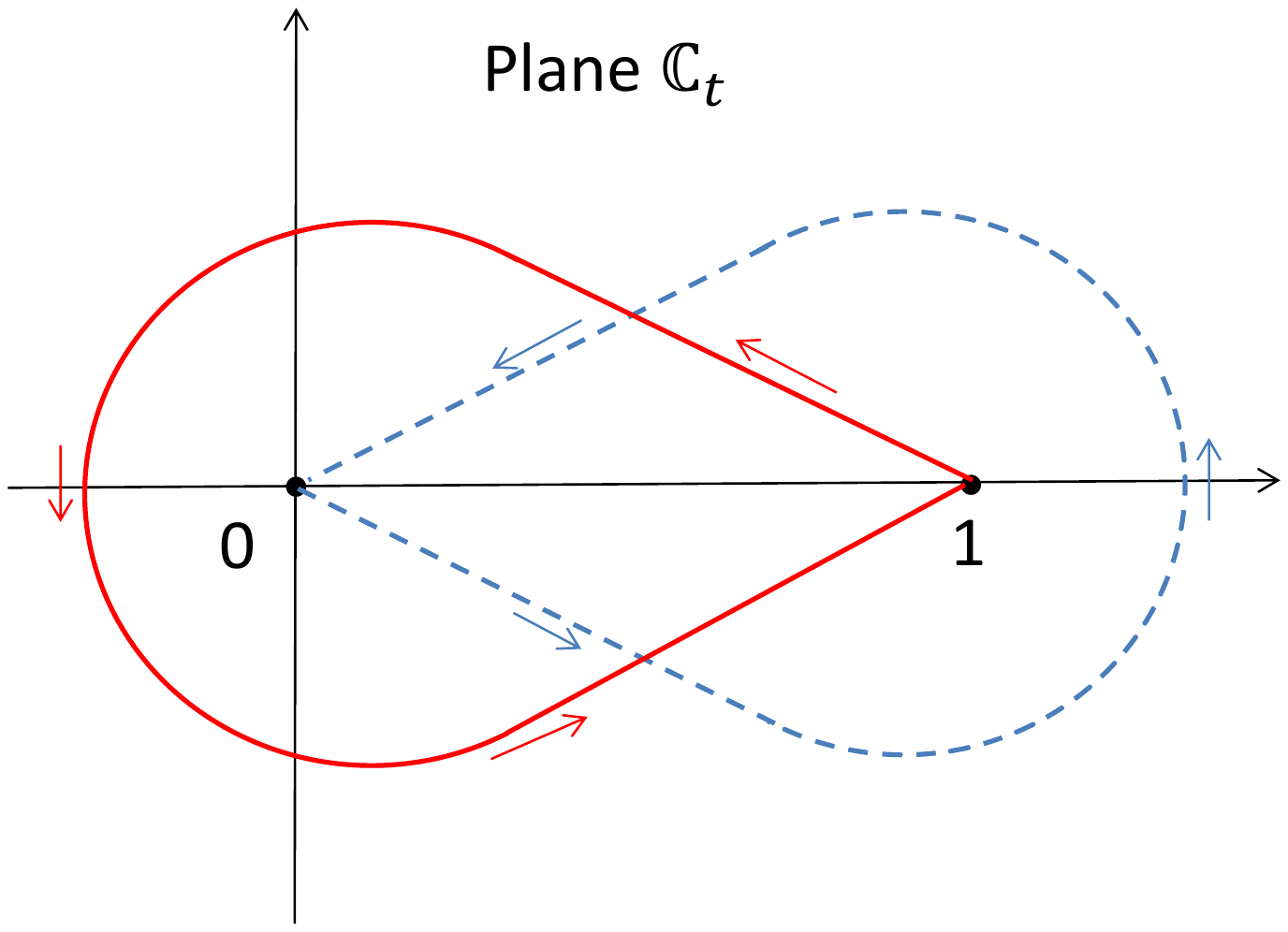}}
\caption{\textit{Integration contours around points 0 and 1.}}
\label{IC}
\end{figure}

\begin{theorem}
\textbf{Let $\nu \in \mathbb{C}$ with $\mathrm{Re}(\nu) < 0$. Then}
\begin{itemize}
\item[\textbf{a)}] \textbf{the operator $\mathfrak{L}:\mathscr{H}_0 \rightarrow \mathscr{H}_0$ is a bijection;}
\item[\textbf{b)}] \textbf{given $K \in \mathscr{H}_0$, the unique solution $E^* = \mathfrak{L}^{-1}K \in \mathscr{H}_0$ to the integral equation (\ref{EI0}) has the integral representation}
\begin{align}
E^*(z) & = \, \mathfrak{L}^{-1} \, K(z) 
\nonumber \\
& = \, \frac{1-x}{2i\pi x} \, e^{z} 
\int_1^{(0+)} \frac{e^{-xtz}}{t(t-1)} \, 
K \left (z \, (-t)^{\nu}(1-t)^{1-\nu} \right ) \, \mathrm{d}t,
\quad z \in \mathbb{C},
\label{SolEI0bis}
\end{align}
\end{itemize}
\textbf{where the contour in integral (\ref{SolEI0}) in variable $t$ is a 
loop starting and ending at point $t = 1$, and encircling the origin 
$t = 0$ once in the positive sense (see Fig.\ref{IC}, red solid line).}
\label{C4}
\end{theorem}

\begin{proof} 
\textbf{a)} Given $K \in \mathscr{H}_0$, equation (\ref{EI0}) for $E^* \in \mathscr{H}_0$ is equivalent to system (\ref{T0}) for the coefficients $(E_\ell)_{\ell \geqslant 1}$ of the exponential series expansion of $E^*$. For $\mathrm{Re}(\nu) < 0$, Corollary \ref{C2} entails these coefficients are uniquely determined by expression (\ref{E0}). The linear operator 
$\mathfrak{L}:\mathscr{H}_0 \rightarrow \mathscr{H}_0$ is consequently one-to-one and onto, and has an inverse $\mathfrak{L}^{-1}$ on $\mathscr{H}_0$.

\textbf{b)} An integral representation for the inverse operator 
$\mathfrak{L}^{-1}$ is now derived as follows. Setting $S = E$ and 
$T = \widetilde{K}$ in (\ref{GsExp}), with the sequence 
$\widetilde{K} = (\widetilde{K}_b)_{b \geqslant 1}$ defined as in 
(\ref{S1}), we obtain
\begin{align}
\mathfrak{G}_{E}^*(z) & = e^z \cdot 
\sum_{b \geqslant 1} (-1)^b \widetilde{K}_b \, \frac{z^b}{b!} \, 
\Phi(b\nu;b;-x \, z)
\nonumber \\
& = \frac{x-1}{x} \, e^z \cdot 
\sum_{b \geqslant 1} (-1)^b \frac{\Gamma(b - b \nu)}
{\Gamma(b)\Gamma(1-b \nu)} \cdot K_b \, \frac{z^b}{b!} \, 
\Phi(b\nu;b;-x \, z)
\label{Sol0}
\end{align}
after using (\ref{S1}) to express $\widetilde{K}_b$ in terms of $K_b$, 
$b \geqslant 1$. Invoke then the integral representation
\begin{equation}
\Phi(\alpha;\beta;Z) = - \frac{1}{2 i \pi}
\frac{\Gamma(1-\alpha)\Gamma(\beta)}{\Gamma(\beta-\alpha)} 
\int_1^{(0+)} e^{Zt}(-t)^{\alpha - 1}(1-t)^{\beta-\alpha-1} \, 
\mathrm{d}t
\label{IRK1bis}
\end{equation}
of the Confluent Hypergeometric function $\Phi(\alpha;\beta;\cdot)$ for 
$\mathrm{Re}(\beta - \alpha) > 0$ \cite[Sect.6.11.1, (3)]{ERD81}, where the integration contour is specified as in Fig.\ref{IC}, red solid line. On account of (\ref{IRK1bis}) applied to $\alpha = b\nu$ and 
$\beta = b \in \mathbb{N}^*$ with $\mathrm{Re}(\nu) < 1$, expression 
(\ref{Sol0}) now reads
\begin{align}
\mathfrak{G}^*_{E}(z) & = \frac{1-x}{2i\pi \, x} \, e^{z} 
\int_1^{(0)^+} \frac{e^{-xzt} \, \mathrm{d}t}{t(t-1)} 
\sum_{b \geqslant 1} (-1)^b \frac{K_b}{b!} \, 
(z(-t)^{\nu}(1-t)^{1-\nu})^b
\nonumber \\
& = \frac{1-x}{2i\pi x} \, e^{z} 
\int_1^{(0+)} \frac{e^{-xtz}}{t(t-1)} \, 
K \left (z \, (-t)^{\nu}(1-t)^{1-\nu} \right ) \, \mathrm{d}t
\label{Sol1bis}
\end{align}
for all $z \in \mathbb{C}$. As $E^*(z) = \mathfrak{G}_E^*(z)$ by definition, expression (\ref{Sol1bis}) readily yields the final representation (\ref{SolEI0bis}), as claimed.
\end{proof}

\noindent
As mentioned in the latter proof, the representation (\ref{SolEI0bis}) of the inverse $\mathfrak{L}^{-1}$ is actually valid for $\mathrm{Re}(\nu) < 1$, although the operator $\mathfrak{L}$ is defined on space 
$\mathscr{H}_0$ for $\mathrm{Re}(\nu) < 0$ only. 

Now, using suitable variable changes in formula (\ref{SolEI0bis}), alternative integral representations for $\mathfrak{L}^{-1}$ with $\mathrm{Re}(\nu) < 0$ can be asserted as follows.

\begin{figure}[htb]
\scalebox{0.8}{\includegraphics[width=14cm, trim = 1.5cm 21cm 3cm 
2.5cm,clip]
{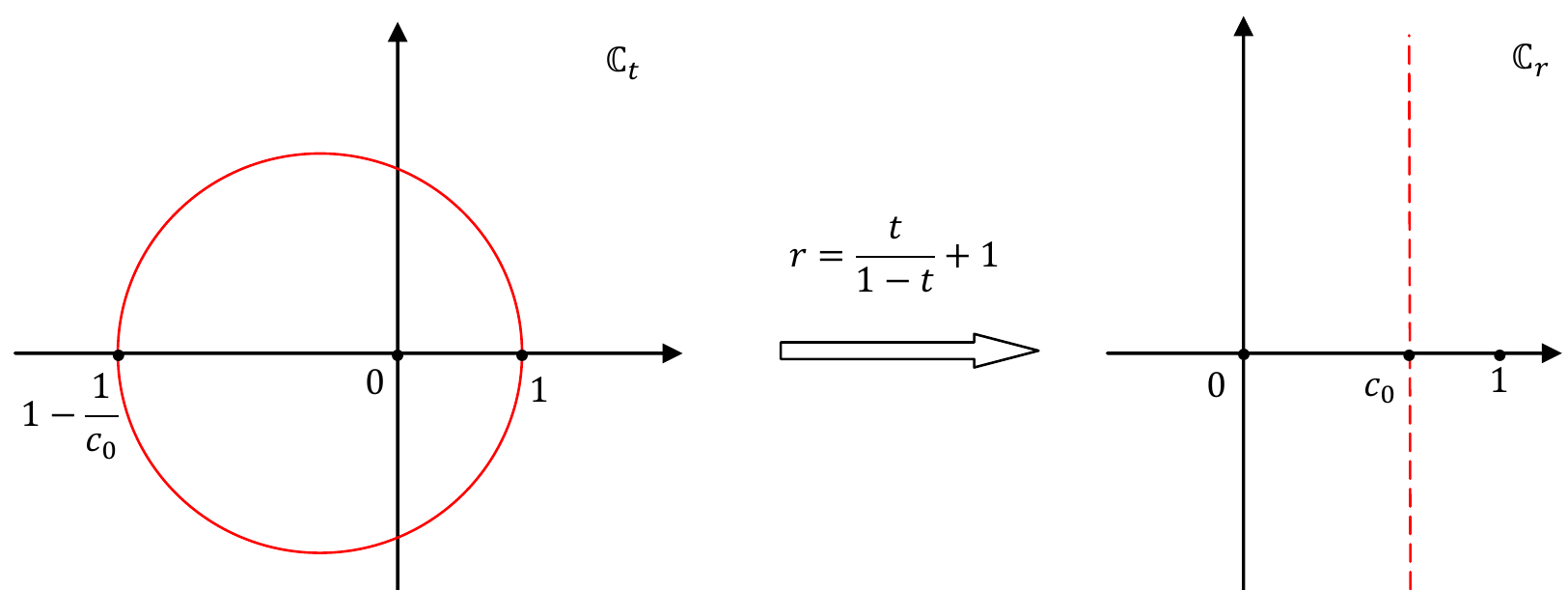}}
\caption{\textit{Transformed integration contours.}}
\label{MAP}
\end{figure}

\begin{corol}
\textbf{Let $\nu \in \mathbb{C}$ with $\mathrm{Re}(\nu) < 0$.}

\textbf{Given $K \in \mathscr{H}_0$, the unique solution 
$E^* = \mathfrak{L}^{-1}K \in \mathscr{H}_0$ to integral equation 
(\ref{EI0}) has the equivalent integral representations}
\begin{align}
E^*(z) & = \, \mathfrak{L}^{-1} \, K(z) 
\nonumber \\
& = \, \frac{1-x}{2i\pi x} \, e^{(1-x)z} 
\int_0^{(1+)} \frac{e^{xtz}}{t(1-t)} \, 
K \left (z \, t^{1-\nu}(t-1)^\nu \right ) \, \mathrm{d}t,
\quad z \in \mathbb{C},
\label{SolEI0}
\end{align}
\textbf{where the contour in (\ref{SolEI0}) in variable $t$ is a loop starting and ending at point $t = 0$, and encircling point $t = 1$ once in the positive sense (see Fig.\ref{IC}, blue dotted line), and}
\begin{align}
E^*(z) & = \, \mathfrak{L}^{-1} \, K(z) 
\nonumber \\
& = \, \frac{1-x}{2i\pi \, x} \;  
e^{(1-x)z} \int_{c_0-i\infty}^{c_0+i\infty} 
\frac{e^{\frac{x \, z}{r}}}{1-r} 
\; K \left (\frac{z(1-r)^\nu}{r} \right ) \, \mathrm{d}r,
\quad z \in \mathbb{C},
\label{SolEI0bisV}
\end{align}
\textbf{where the contour in (\ref{SolEI0bisV}) is the vertical line $\mathrm{Re}(r) = c_0$, for any real abcissa $0 < c_0 < 1$ (see Fig.\ref{MAP}, red dotted line).}
\label{C4bis}
\end{corol}

\begin{proof}
$\bullet$ By the variable change $t \mapsto 1-t$, formula  
(\ref{SolEI0bis}) readily entails (\ref{SolEI0}) which is defined for 
$\mathrm{Re}(\nu) < 0$.

$\bullet$ Let $0 < c_0 < 1$. As a contour in integral (\ref{SolEI0bis}), choose the circle centered at point $1 - 1/(2c_0)$ on the real axis and with radius $1/(2c_0)$; this circle passes through point $1$ and encircles the origin (see Fig.\ref{MAP}). It is easily verified that the homographic transformation $t \mapsto r' = t/(1-t)$ maps this circle 1-to-1 and onto the vertical line $\mathrm{Re}(r') = c_0 -1$.  Applying the latter variable change $t \mapsto r'$ to (\ref{SolEI0bis}) with 
$$
(-t)^{\nu}(1-t)^{1-\nu} = \frac{(-r')^\nu}{1+r'}, \qquad  
\mathrm{d}t = \frac{\mathrm{d}r'}{(1+r')^2}
$$
then readily gives 
$$
\mathfrak{L}^{-1}K(z) = \left ( \frac{1-x}{2i\pi x} \right) \, 
e^{z(1-x)} \int_{\mathrm{Re}(r') = c_0-1} e^{\frac{z \, x}{1+r'}} \, 
K \left ( \frac{z (-r')^\nu}{1+r'} \right ) \frac{\mathrm{d}r'}{r'}
$$
for all $z \in \mathbb{C}$. The mapping $r' \mapsto r = r' + 1$ then eventually transforms the latter integral to the expected representation 
(\ref{SolEI0bisV}) with integration contour the vertical line 
$\mathrm{Re}(r) = c_0$, $0 < c_0 < 1$.
\end{proof}

By the factorization (\ref{DefTTbis}), we readily deduce that the inverse of operator $\mathfrak{M}$ is given by
\begin{equation}
\mathfrak{M}^{-1}g(z) = \frac{x \, \nu(1-\nu)}{1-x} \cdot 
\mathfrak{L}^{-1}(z \, g')(z), \qquad z \in \mathbb{C},
\label{InvM}
\end{equation}
for all $g \in \mathscr{H}_0$, with inverse $\mathfrak{L}^{-1}$ provided by either integral representation (\ref{SolEI0bis}), (\ref{SolEI0}) or (\ref{SolEI0bisV}). The involvement of the derivative $g'$ for the inverse $\mathfrak{M}^{-1}g$ in (\ref{InvM}) reminds us of formula (\ref{SolAbel}) in the particular case of the Abel's equation.




\section{Appendix}


\subsection{Proof of Proposition \ref{LemT0} (continued)}
\label{A4}
We conclude the proof of Proposition \ref{LemT0} by expressing the coefficients $Q_{b,\ell}$, $1 \leqslant \ell \leqslant b$, introduced in 
(\ref{DefQbell}) in terms of Hypergeometric polynomials only. We first calculate coefficients $Q_{b,\ell}(s)$, $1 \leqslant \ell \leqslant b$, in terms of the general Gauss Hypergeometric function $F$. Recall that 
$F = F(\alpha,\beta;\gamma;\cdot)$ has the integral representation 
\cite[Chap.15, Sect.15.6.1]{NIST10}
\begin{equation}
F(\alpha,\beta;\gamma;z) = 
\frac{\Gamma(\gamma)}{\Gamma(\beta)\Gamma(\gamma-\beta)}
\int_0^1 \frac{t^{\beta-1}(1-t)^{\gamma-\beta-1}}{(1-z t)^\alpha} 
\, \mathrm{d}t, \quad \vert z \vert < 1,
\label{HyperGauss}
\end{equation}
for real parameters $\alpha$, $\beta$, $\gamma$ where 
$\gamma > \beta > 0$.

\begin{lemma}
\textbf{We have}
\begin{align}
Q_{b,\ell} = & \, 
- \frac{\Gamma(\ell)\Gamma(1 - b \, \nu)}{\Gamma(\ell + 1 - b \, \nu)} 
\, \left ( \frac{x^{1+b}}{1-x}  \right ) \, 
\Bigl [  \, \nu \, (b-\ell) \; \times 
\nonumber \\
& \,F(b \, (1-\nu),\ell;\ell + 1 - b \, \nu;1-x)  + (\ell - b \, \nu) \times 
F(b \, (1-\nu),\ell;\ell - b \, \nu;1-x) 
\, \Bigr ]
\label{CoeffQ}
\end{align}
\textbf{for $1 \leqslant \ell \leqslant b$.}
\label{CoeffT0}
\end{lemma}

\begin{proof}
To calculate the integral $M_{b,\ell}$ introduced in (\ref{DefMbell}), use the definition of $\mathfrak{R}(t)$, to write
$$
M_{b,\ell} = x^{b}\,\int_0^1 
t^\ell (1-t)^{-b \, \nu} 
\left ( 1 - (1-x)t \right )^{b \, (\nu - 1)} \, 
\mathrm{d}t;
$$
using representation (\ref{HyperGauss}) for parameters 
$\alpha = -b(1-\nu)$, $\beta = \ell + 1$, $\gamma = 2 + \ell - b \, \nu$, this integral reduces to
\begin{equation}
M_{b,\ell} = 
\frac{\Gamma(\ell+1)\Gamma(1- b \, \nu)}{\Gamma(2+\ell -  b \, \nu)} \,x^{b}\, 
F(b \, (1-\nu),\ell+1;2+\ell-b \, \nu;1-x);
\label{MBL0}
\end{equation}
after (\ref{MBL0}) and the expression (\ref{DefQbell}) of coefficient 
$Q_{b,\ell}$, we then derive
\begin{align}
Q_{b,\ell} = & \, 
\frac{\Gamma(\ell)\Gamma(1 - b \, \nu)}{\Gamma(\ell + 1 - b \, \nu)} 
\, x^{b}\,\times 
\nonumber \\
& \, \Bigl [  \frac{\ell}{\ell+1-b \, \nu} (\ell+1-b) \cdot 
F(b \, (1-\nu),\ell+1;\ell + 2 - b \, \nu;1-x) 
\nonumber \\
& \; - \ell \, c \cdot F(b \, (1-\nu),\ell;\ell + 1 - b \, \nu;1-x) 
\, \Bigr ].
\label{MBL1}
\end{align}
To simplify further the latter expression, first invoke the identity 
\begin{equation}
\beta \, F(\alpha,\beta+1;\gamma+1;z) = \gamma \, 
F(\alpha,\beta;\gamma;z) - 
(\gamma-\beta) \, F(\alpha,\beta;\gamma+1;z)
\label{ID1}
\end{equation}
easily derived from representation (\ref{HyperGauss}) for 
$F(\alpha,\beta+1;\gamma+1;z)$, after splitting the factor $t^\beta$ of the integrand into $t^\beta = t^{\beta-1} - t^{\beta-1}(1-t)$. Applying 
(\ref{ID1}) to $\alpha = b \, (1-\nu)$, $\beta = \ell$ and 
$\gamma = \ell + 1 - b \, \nu$ then enables one to express the term 
$F(b \, (1-\nu),\ell+1;\ell + 2 - b \, \nu;1-x)$ in the r.h.s. of 
(\ref{MBL1}) as a combination of 
$F(b \, (1-\nu),\ell;\ell + 1 - b \, \nu;1-x)$ and 
$F(b \, (1-\nu),\ell;\ell + 2 - b \, \nu;1-x)$ hence, after simple algebra,
\begin{align}
Q_{b,\ell} = & \, 
\frac{\Gamma(\ell)\Gamma(1 - b \, \nu)}{\Gamma(\ell + 1 - b \, \nu)} 
\, x^{b}\, \Bigl [ \left\{ (\ell + 1 - b) - \ell \, c \right \} \cdot 
F(b \, (1-\nu),\ell;\ell + 1 - b \,\nu;1-x) 
\nonumber \\
& \, - \frac{(\ell+1-b)(1-b \, \nu)}{\ell + 1 - b \, \nu} \cdot 
F(b \, (1-\nu),\ell;\ell + 2 - b \, \nu;1-x) \, \Bigr ].
\label{MBL2}
\end{align}
Furthermore, the contiguity identity \cite[Sect.15.5.18]{NIST10}
\begin{align}
& \gamma[\gamma-1-(2\gamma-\alpha-\beta-1)z] \, F(\alpha,\beta;\gamma;z) 
\; + 
\nonumber \\
& (\gamma-\alpha)(\gamma-\beta)z \, F(\alpha,\beta;\gamma+1;z) = 
\gamma(\gamma-1)(1-z) \, F(\alpha,\beta;\gamma-1;z)
\end{align}
applied to $\alpha = b \, (1-\nu)$, $\beta = \ell$ and 
$\gamma = \ell + 1 - b \, \nu$ allows us to write the last term 
$F(b \, (1-\nu),\ell;\ell + 2 - b \, \nu;1-x)$ in the bracket of the r.h.s. of (\ref{MBL2}) as a combination of $F(b \, (1-\nu),\ell;\ell - b \, \nu;1-x)$ and $F(b \, (1-\nu),\ell;\ell + 1 - b \, \nu;1-x)$, that is,
\begin{align}
& F(b \, (1-\nu),\ell;\ell + 2 - b \, \nu;1-x) = 
\frac{\ell + 1 - b \, \nu}{(\ell + 1 - b)(1-b \, \nu)(1-x)} \; \times 
\nonumber \\
& \Bigl [ (\ell - b \, \nu)x \cdot 
F(b \, (1-\nu),\ell;\ell - b \, \nu;1-x) \; - 
\nonumber \\
& [\ell - b \, \nu \; - (\ell + 1 - b - b \, \nu)(1-x)] \cdot 
F(b \, (1-\nu),\ell;\ell + 1 - b \, \nu;1-x) \Bigr ];
\nonumber
\end{align}
inserting the latter relation into the right-hand side of (\ref{MBL2}) then yields
\begin{align}
Q_{b,\ell} = & \, 
\frac{\Gamma(\ell)\Gamma(1 - b \, \nu)}{\Gamma(\ell + 1 - b \, \nu)} 
\, x^{b}\,\times \Bigl [ T_{b,\ell} \cdot 
F(b \, (1-\nu),\ell;\ell + 1 - b \, \nu;1-x) \; -  
\nonumber \\
& \, \frac{1}{1-x}(\ell - b \, \nu)x \cdot 
F(b \, (1-\nu),\ell;\ell - b \, \nu;1-x) \, \Bigr ]
\label{MBL3}
\end{align}
where 
$$
T_{b,\ell} = b \, \nu - \ell \, c + \frac{(\ell - b \, \nu)}{1-x} = 
\frac{(\ell - b) \nu x}{1-x}
$$
after the definition $c = (1-\nu x)/(1-x)$ of constant $c$. Inserting this value of $T_{b,\ell}$ in the right-hand side of (\ref{MBL3}) readily provides expression (\ref{CoeffQ}) for $Q_{b,\ell}$. 
\end{proof}

We finally show how coefficient $Q_{b,\ell}$ can be written in terms of a Hypergeometric polynomial only. Applying the general identity 
\cite[Chap.9, Sect.9.131.1]{GRAD07}
\begin{equation}
F(\alpha,\beta;\gamma;z) = (1-z)^{\gamma-\alpha-\beta}
F(\gamma-\alpha,\gamma-\beta;\gamma;z), \qquad \vert z \vert < 1,
\label{GI0}
\end{equation}
to each term $F(b(1-\nu),\ell;\ell+1-b\nu;1-x)$ and 
$F(b(1-\nu),\ell;\ell-b\nu;1-x)$ in (\ref{CoeffQ}), we obtain
\begin{equation}
Q_{b,\ell} = 
- \frac{\Gamma(\ell)\Gamma(1 - b \, \nu)}{\Gamma(\ell + 1 - b \, \nu)} 
\cdot \frac{x^2}{b (1-x)} (\ell - b\nu) \, \times 
R_{b,\ell}, 
\qquad b \geqslant \ell \geqslant 1, 
\label{Q1}
\end{equation}
where we set
\begin{align} 
R_{b,\ell} = & \; \frac{b\nu \, (b-\ell)}{\ell - b \nu} \cdot 
F(\ell - b + 1,-b\nu + 1;\ell - b\nu + 1;1-x) \; + 
\nonumber \\
& \; b \, x^{-1} \cdot F(\ell-b,-b\nu;\ell-b\nu;1-x).
\nonumber
\end{align}
From the identity \cite[Chap.15, Sect.15.5.1]{NIST10}
\begin{equation}
\frac{\mathrm{d}}{\mathrm{d}z} \, F(\alpha,\beta;\gamma;z) = 
\frac{\alpha \beta}{\gamma} \, 
F(\alpha+1,\beta+1;\gamma+1;z), \qquad \vert z \vert < 1,
\label{DF}
\end{equation}
applied to parameters $\alpha = \ell - b$, $\beta = -b\nu$ and 
$\gamma= \ell - b\nu$, the factor $R_{b,\ell}$ above then equals the derivative 
\begin{align}
R_{b,\ell} = & \; x^{b}\,\frac{\mathrm{d}}{\mathrm{d}z} (1-z)^{-b}F(\ell-b,-b\nu;\ell-b\nu;z) \vert_{z = 1-x} = 
x^{b}\,\frac{\mathrm{d}}{\mathrm{d}z} F(b-b\nu,\ell;\ell-b\nu;z) \vert_{z = 1-x}
\nonumber \\
= & \; x^{b}\,\frac{(b-b\nu)\ell}{\ell-b\nu} \, F(b-b\nu+1,\ell+1;\ell-b\nu+1;1-x)
\nonumber
\end{align}
hence
\begin{equation}
R_{b,\ell} = \frac{(b-b\nu)\ell}{\ell-b\nu} \, x^{-1}
F(\ell - b,-b\nu;\ell-b\nu+1;1-x)
\label{Q1bis}
\end{equation}
where we have successively applied identity (\ref{GI0}), (\ref{DF}) and 
(\ref{GI0}) again to derive the second, third and fourth equality, respectively. Using (\ref{Q1bis}), expression (\ref{Q1}) for 
$Q_{b,\ell}$ then reads
\begin{equation}
Q_{b,\ell} = - \frac{\Gamma(\ell)\Gamma(1-b\nu)}{\Gamma(\ell-b\nu)} \, 
\frac{x}{1-x} \; 
\frac{\ell(1-\nu)}{\ell - b\nu} \; S_{b,\ell}(\nu;1-x), 
\; \; b \geqslant \ell \geqslant 1,
\label{Q3}
\end{equation}
where we set
$$
S_{b,\ell}(\nu;1-x) = F(\ell-b,-b\nu;\ell-b\nu+1;1-x).
$$
To reduce further $S_{b,\ell}(\nu;1-x)$, invoke the identity 
\cite[Chap.15, Sect.15.8.7]{NIST10}
\begin{equation}
F(-m,\beta,\gamma;1-x) = 
\frac{\Gamma(\gamma)\Gamma(\gamma-\beta+m)}
{\Gamma(\gamma-\beta)\Gamma(\gamma+m)} \, F(-m,\beta,\beta+1-m-\gamma;x), 
\quad x \in \mathbb{C},
\label{IDFG}
\end{equation}
for any non negative integer $m$ and complex numbers $\beta$, $\gamma$ such that $\mathrm{Re}(\gamma) > \mathrm{Re}(\beta)$; applying 
(\ref{IDFG}) to factor $S_{b,\ell}(\nu;1-x)$ in (\ref{Q3}) then readily gives the final expression (\ref{Q0}) for all indexes 
$b \geqslant \ell \geqslant 1$. This concludes the proof of Proposition \ref{LemT0} $\blacksquare$

\subsection{Proof of Corollary \ref{C3}}
\label{A5}
$\bullet$ From the definition (\ref{DefSS}) of integral operator 
$\mathfrak{M}$, split the integral 
\begin{align}
\mathfrak{M}E^*(z) = & \, 
\int_0^1 e^{- z \, t^{-\nu}(1-(1-x)t)} \, 
E^* \left (z t^{-\nu}(1-t) \right ) \frac{\mathrm{d}t}{t} 
\nonumber \\
= & \, \int_0^{\widehat{t}} \; (...) \, \frac{\mathrm{d}t}{t} + 
\int_{\widehat{t}}^1 \; (...) \, \frac{\mathrm{d}t}{t}
\nonumber
\end{align}
over adjacent segments $[0,\widehat{t}]$ and $[\widehat{t},1]$, respectively; applying the variable change $\tau = t^{-\nu}(1-t)$ on each of these two intervals with 
$\tau = \tau_-(t) \Leftrightarrow t = 
\theta_-(\tau) \in [0,\widehat{t}]$ and 
$\tau = \tau_+(t) \Leftrightarrow t = 
\theta_+(\tau) \in [\widehat{t},1]$ by the definition (\ref{Thetapm}) of mappings $\theta_-$ and $\theta_+$, we then successively obtain
\begin{align}
\mathfrak{M}E^*(z) = & \, \int_0^{\widehat{\tau}} 
e^{-z\theta_-(\tau)^{-\nu}(1-(1-x)\theta_-(\tau))} 
E^* \left (z \, \tau \right ) \, 
\frac{-\mathrm{d}\tau}{\theta_-(\tau)^{-\nu}(\nu+(1-\nu)\theta_-(\tau))} 
\nonumber \\
+ & \, \int_{\widehat{\tau}}^0
e^{-z \theta_+(\tau)^{-\nu}(1-(1-x)\theta_+(\tau))} 
E^* \left (z \, \tau \right ) \, 
\frac{-\mathrm{d}\tau}{\theta_+(\tau)^{-\nu}(\nu+(1-\nu)\theta_+(\tau))}
\nonumber
\end{align}
with $\widehat{\tau} = \tau_-(\widehat{t}) = \tau_+(\widehat{t})$ and the differential 
$\mathrm{d}t/t = -\mathrm{d}\tau/[t^{-\nu}(\nu+(1-\nu)t)]$; this readily reduces to a single integral over segment $[0,\widehat{\tau}]$, that is,
$$
\mathfrak{M}E^*(z) = \int_0^{\widehat{\tau}} 
\left [ \Psi_-(z,\tau) - \Psi_+(z,\tau) \right ]  
E^* \left ( z \, \tau \right ) \, \mathrm{d}\tau 
$$
with $\Psi_-(z,\tau)$ and $\Psi_+(z,\tau)$ given as in the Corollary. The final variable change $\xi = z \cdot \tau$ yields the right-hand side  of (\ref{VOL0}) and the corresponding integral equation. 

$\bullet$ We finally verify that the r.h.s. of (\ref{VOL0}) is well-defined for any $E^* \in \mathscr{H}_0$. The denominator 
$t^{-\nu}(-\nu + (\nu-1)t)$ of $\Psi_-(z,\tau)$ with 
$t = \theta_-(\tau)$ (resp. of $\Psi_+(z,\tau)$ with 
$t = \theta_+(\tau)$) vanishes at either $\tau = 0$ or 
$\tau = \widehat{\tau}$ (resp. at $\tau = \widehat{\tau}$). As to the possible singularity at $\tau = 0$ for $\Psi_-(z,\tau)$, we have 
$\tau \sim t^{-\nu}$ for small $t = \theta_-(\tau)$ so that 
$$
\frac{1}{t^{-\nu}(-\nu + (\nu-1)t)} \sim - \, \frac{t^{\nu}}{\nu} \sim 
- \frac{1}{\nu\tau}, \qquad \tau \rightarrow 0;
$$
the product $E^*(z\tau) \cdot \Psi_-(z,\tau)$ is thus integrable near 
$\tau = 0$ for any $E^* \in \mathscr{H}_0$. Furthermore, a second order Taylor expansion of $\tau = \tau(t)$ near point $t = \widehat{t}$ yields 
$\tau = \widehat{\tau} + \tau''(\widehat{t}) \, (t-\widehat{t})^2/2 
+ o(t-\widehat{t})^2$ 
with $\tau'(\widehat{t}) = 0$ by definition and 
$\tau''(\widehat{t}) < 0$; as a result,
$$
t - \widehat{t} \sim 
\pm \sqrt{\frac{-2(\widehat{\tau}-\tau)}{\tau''(\widehat{t})}}, 
\qquad \tau \rightarrow \widehat{\tau}.
$$
The denominator $t^{-\nu}(-\nu+(\nu-1)t)$ of either $\Psi_-(z,\tau)$ or 
$\Psi_+(z,\tau)$ is consequently asymptotic to
$$
t^{-\nu}(-\nu+(\nu-1)t) \sim 
(\widehat{t})^{-\nu}(\nu-1)(t - \widehat{t}) \sim 
\pm (\widehat{t})^{-\nu}(\nu-1) 
\sqrt{\frac{2(\widehat{\tau}-\tau)}{-\tau''(\widehat{t})}}
$$
when $\tau \rightarrow \widehat{\tau}$; the singularity of 
$\Psi_-(z,\tau)$ (resp. $\Psi_+(z,\tau)$) at point $\tau = \widehat{\tau}$ is consequently of order 
$$
\Psi_-(z,\tau) = O \left ( \frac{1}{\sqrt{\widehat{\tau}-\tau}} \right ), 
\quad 
\Psi_+(z,\tau) = O \left ( \frac{1}{\sqrt{\widehat{\tau}-\tau}} \right )
$$
and the kernel $\Psi(z,\cdot) = \Psi_-'(z,\cdot) - \Psi_+(z,\cdot)$ is thus integrable at $\tau = \widehat{\tau}$. This ensures that the singular integral (\ref{VOL0}) is well-defined for any 
$E^* \in \mathscr{H}_0$ $\blacksquare$


\end{document}